\documentclass[11pt]{article}%
\usepackage{stix2}
\usepackage{amsfonts}
\usepackage{amsmath}
\usepackage{amssymb}
\usepackage{amsxtra,amscd}
\usepackage[amsmath,hyperref,thmmarks]{ntheorem}
\usepackage{makeidx}
\usepackage{graphicx}
\usepackage[colorlinks=true,bookmarksnumbered=true,pdfpagemode=None,final]%
{hyperref}%
\providecommand{\U}[1]{\protect\rule{.1in}{.1in}}
\theoremnumbering{arabic}
\theoremheaderfont{\scshape}
\RequirePackage{latexsym}
\theorembodyfont{\slshape}
\theoremseparator{.}
\newtheorem*{conjecture}{Conjecture}
\newtheorem*{theorem}{Theorem}
\theorembodyfont{\upshape}
\theorembodyfont{\small}
\newtheorem*{aside}{Aside}

\theorembodyfont{\normalsize}
\theoremstyle{nonumberplain}
\theoremsymbol{\ensuremath{_\Box}}

\qedsymbol{\ensuremath{_\Box}}
\newcommand{\df}{\smash{\lower.12em\hbox{\textup{\tiny def}}}}
\newcommand{\lsimeq}{\smash{\lower.24em\hbox{$\scriptstyle\simeq$}}}

\DeclareMathOperator{\Gal}{Gal}
\DeclareMathOperator{\Spec}{Spec}
\DeclareMathOperator{\Hom}{Hom}
\DeclareMathOperator{\Ext}{Ext}
\DeclareMathOperator{\Ker}{Ker}
\DeclareMathOperator{\Br}{Br}
\DeclareMathOperator{\Pic}{Pic}
\newcommand{\dstyle}{\displaystyle}
\newcommand{\textnf}{\textnormal}
\def\1{\mathbb{1}}

\DeclareFontEncoding{OT2}{}{} 
\newcommand{\textcyr}[1]{%
{\fontencoding{OT2}\fontfamily{cmr}\fontseries{m}\fontshape{n}\selectfont #1}}
\newcommand{\TS}{{\mbox{\textcyr{Sh}}}}

\def\mapsto{\DOTSB\mathchar"39AD }

\usepackage{changepage}

\usepackage[overload]{textcase}

\usepackage{xcolor}
\definecolor{bblue}{rgb}{0.0, 0.0, 0.6}

\usepackage{array}

\usepackage{microtype}

\usepackage{tikz}
\usepackage{tikz-cd}
\usetikzlibrary{matrix,arrows,positioning,decorations.pathmorphing}
\usetikzlibrary{decorations.markings,shapes.geometric}
\tikzset{commutative diagrams/column sep/Huge/.initial=12ex}
\tikzcdset{arrow style=math font}

\usepackage[shortlabels]{enumitem}
\setitemize[1]{label=$\diamond$\hspace{0.07in}}
\setlist{topsep=0.1em,itemsep=0.1em,parsep=0.1em}

\setdescription{font=\normalfont}

\usepackage{url}
\usepackage{verbatim}

\usepackage[notcite,notref]{showkeys}

\usepackage{mathtools}

\newcommand{\eb}[1]{{\itshape\bfseries#1}}
\renewcommand{\emph}{\eb}

\usepackage{titletoc,titlesec}
\contentsmargin{2.0em}
\dottedcontents{section}[2.3em]{}{1.8em}{1pc}
\usepackage[nottoc]{tocbibind}

\titleformat*{\section}{\LARGE\bfseries}
\titleformat*{\subsection}{\Large\itshape}
\titleformat*{\subsubsection}{\scshape}

\hypersetup{linkcolor=bblue,anchorcolor=bblue,citecolor=bblue,filecolor=bblue,urlcolor=bblue}
\hypersetup{pdftitle={Arithmetic Duality},pdfauthor={J.S. Milne}}
\hypersetup{pdfsubject={none},pdfkeywords={none}}

\usepackage[papersize={6.5in,10.0in},margin=0.5in,pdftex]{geometry}

\begin{document}

\title{Arithmetic Duality}
\author{J.S.\ Milne}
\date{\today}
\maketitle

\begin{abstract}
In the 1950s and 1960s Tate proved some duality theorems in the Galois
cohomology of finite modules and abelian varieties. As for most of Tate's work
this has had a profound influence on mathematics with many applications and
further developments. In this article, I discuss Tate's theorems and some of
these developments.

\end{abstract}
\tableofcontents

\renewcommand{\thefootnote}{\fnsymbol{footnote}} \footnotetext{This is the
written version of my talk at the conference, \textit{The Legacy of John Tate, and
Beyond}. Harvard University, March 17--21, 2025. I thank the organizers for
inviting me to speak on the ``theme of arithmetic duality.'' Also,
I thank the referee for numerous helpful comments.} \renewcommand{\thefootnote}{\arabic{footnote}}

\bigskip
Since the work of Tate and others in the late 1950s, duality theorems
have played a vital role in arithmetic geometry, similar to that 
of Poincar\'e duality in other areas of mathematics.
As we shall see, mathematicians have often found, at a crucial point in their 
research, that they need such a theorem. Here we describe
the early years of arithmetic duality.

\bigskip\noindent\textsl{Notation:} I generally follow the notation I learned
from Tate. For example, $X(\ell)$ is the $\ell$-primary component of an
abelian group $X$, $X_{n}=\{x\in X\mid nx=0\}$, and $[X]$ is the 
order of $X$. For a field $K$, $\bar K$ denotes a separable closure of $K$.
All varieties are connected.

\section{Local duality (Tate 1957)}

For an abelian variety $A$ over a field $K$, the Galois cohomology group
$H^{1}(K,A)$ classifies the principal homogeneous spaces (torsors) of $A$ over
$K$. Ch\^{a}telet demonstrated the importance of this group in the diophantine
study of elliptic curves and Weil for a general abelian variety, and so, in
his 1957 Bourbaki seminar, Tate named it the Weil--Ch\^{a}telet group.

Except that they are torsion, almost nothing was known about the groups until
Tate proved in his talk that, when $K$ is a local field of characteristic zero
(i.e., a finite extension of $\mathbb{Q}{}_{p}$), the Weil--Ch\^{a}telet group
of $A$ (with its discrete topology) is dual to the group of rational points on
the dual abelian variety $A^{\prime}$ (a compact group), i.e.,%
\[
H^{1}(K,A)\xrightarrow{\lsimeq}A^{\prime}(K)^{\ast},\quad\ast
\,=\text{Pontryagin dual.}%
\]

Since $A^{\prime}(K)$ was well-understood at the time, for example, it
contains a subgroup of finite index isomorphic to $\mathcal{O}_{K}^{\dim A}$
(Lutz, Mattuck), this tells us a great deal about the Weil--Ch\^{a}telet group.
For example, it shows that the non-$p$ part of $H^{1}(K,A)\frac{{}}{{}}$ is
finite, and has a description in terms of the torsion subgroup of $A^{\prime
}(K)$. Lang and Tate had proved that earlier, and it was by thinking about
this and investigating the elliptic curve case that Tate was led to his theorem.

Many readers will recognize the statement as being part of what we now call
Tate local duality. It took Tate some time to realize this. In the final
paragraph of his Bourbaki talk, almost as an afterthought, he noted that there
are canonical pairings for all $r,s$,
\[
H^{r}(K,A)\times H^{s}(K,A^{\prime})\rightarrow H^{r+s+1}(K,\mathbb{G}_{m}),
\]
which he later showed give a duality,%
\[
H^{r}(K,A)\times H^{s}(K,A^{\prime})\rightarrow H^{2}(K,\mathbb{G}_{m}%
)\simeq\mathbb{Q}/\mathbb{Z}{},
\]
when $r+s=1$.

From the exact sequence
\[
0\rightarrow M\rightarrow A(\bar{K})\rightarrow B(\bar{K})\rightarrow0\
\]
defined by an isogeny $A\rightarrow B$ of abelian varieties and its dual,
\[
0\rightarrow M^{D}\rightarrow B^{\prime}(\bar{K})\rightarrow A^{\prime
}(\bar{K})\rightarrow0,\qquad M^{D}\overset{\df}{=}\Hom(M,\bar{K}^{\times}),
\]
Tate deduced a commutative diagram
\[
\begin{tikzcd}[column sep=scriptsize]
H^{0}(K,A)\arrow{r}\arrow{d}{\simeq}
& H^{0}(K,B)\arrow{r}\arrow{d}{\simeq}
&H^{1}(K,M)\arrow{r}\arrow[dashed]{d}
&H^{1}(K,A)\arrow{r}\arrow{d}{\simeq}
& H^{1}(K,B)\arrow{d}{\simeq}\\
H^{1}(K,A^{\prime})^{\ast}\arrow{r}
&H^{1}(K,B^{\prime})^{\ast}\arrow{r}
& H^{1}(K,M^D)^{\ast}\arrow{r}
& H^{0}(K,A^{\prime})^{\ast}\arrow{r}& H^{0}(K,B^{\prime})^{\ast},\end{tikzcd}
\]
and hence a duality%
\[
H^{r}(K,M)\times H^{2-r}(K,M^{D})\rightarrow H^{2}(K,\mathbb{G}_{m}%
)\simeq\mathbb{Q}/\mathbb{Z}.
\]

For some time, Tate thought that the existence of this duality was a curious
property of the Galois submodules of $A(K)$.\footnote{Personal communication.}
Eventually, of course, he realized that all finite Galois modules for a local
field of characteristic zero have this property, and so obtained Tate local
duality: there are compatible dualities
\begin{align*}
H^{r}(K,M)\times H^{2-r}(K,M^{D})  &  \rightarrow H^{2}(K,\mathbb{G}%
_{m})\simeq\mathbb{Q}/\mathbb{Z}\\
H^{r}(K,A)\times H^{1-r}(K,A^{\prime})  &  \rightarrow H^{2}(K,\mathbb{G}%
_{m})\simeq\mathbb{Q}/\mathbb{Z}%
\end{align*}
for $K$ a finite extension of $\mathbb{Q}_{p}$, $M$ a finite Galois module
with dual $M^{D}=\Hom(M,\bar{K}^{\times})$, and $A$ an abelian variety with
dual $A^{\prime}$.

Thus, the duality theorem for abelian varieties was proved before the
(easier!) duality theorem for finite Galois modules, and even before a local
duality theorem was available for elliptic curves.

\subsubsection{The dual abelian variety}

Every abelian variety has a dual, which is an abelian variety of the same
dimension, but not necessarily isomorphic. The dual of an elliptic curve is
the curve itself. Usually, the dual of $A$ is defined to classify translation
invariant line bundles on $A$, but, as Weil observed, when you remove the
zero-section of such line bundle, it acquires a group structure that makes it
an extension of $A$ by $\mathbb{G}_{m}$. In this way, we get an isomorphism%
\[
A^{\prime}(k)\simeq\Ext_{k}(A,\mathbb{G}_{m})\quad\text{(Barsotti-Weil
formula).}%
\]
This interpretation of $A^{\prime}$ makes it easier to define the pairings.
Indeed, when his collected works were published almost 60 years after he gave
his Bourbaki seminar, Tate added a note saying exactly that:

\begin{quote}
In hindsight, the [cohomological] pairing for dual abelian varieties $A$ and
$B$ is evident from the relation $B=\Ext(A,\mathbb{G}_{m})$ (Tate, Collected
Works, Part I, p.127).
\end{quote}

\section{Global duality (Tate 1962)}

Tate immediately recognized the importance of extending his local duality
theorems to global fields. By 1960 he knew the statements he wanted, but not
the proofs. By early 1962 he had the proofs, in time to announce his theorems
at the 1962 ICM in Stockholm.

One statement of his theorem is that there is a nine-term exact sequence, as
below. To understand the sequence, note that the $\beta$'s map the global
Galois cohomology group into a product of the local groups. One would like to
know the kernels and cokernels of these maps, but there is no simple
expression for these. The best one can do is Tate's exact sequence,%

\[
\begin{tikzcd}[column sep=scriptsize]
0 \arrow{r}& H^{0}(K,M) \arrow{r}{\beta^{0}}&
\prod_vH^{0}(K_v,M)\arrow{r}{\gamma^{0}}
&H^{2}(K,M^{D})^{\ast}\arrow{d}\\
&H^{1}(K,M^{D})^{\ast}\arrow{d}
&\mathclap{\coprod}\mathclap{\prod}\,\,\,{}_v H^{1}(K_v,M)\arrow{l}[swap]{\gamma^{1}}
&H^{1}(K,M)\arrow{l}[swap]{\beta^{1}}\\
&H^{2}(K,M) \arrow{r}{\beta^{2}}
&\bigoplus_vH^{2}(K_v,M)\arrow{r}{\gamma^{2}}
&H^{0}(K,M^{D})^{\ast}\arrow{r}
&0,\end{tikzcd}
\]
where

\begin{itemize}
\item $K$ = global field; $\bar{K}$ = separable closure of $K$;

\item $M$ finite $\Gal(\bar{K}/K)$-module; $\mathrm{char}(K)\nmid\lbrack M]$
(order of $M$);

\item $M^{D}=\Hom(M,\bar{K}^{\times})$; $\quad\ast$ = Pontryagin dual;

\item $H^{0}(\mathbb{R}{},M)=M^{\Gal(\mathbb{C}{}/\mathbb{R}{})}%
/(1+\iota)M,\quad H^{0}(\mathbb{C}{},M)=0$;

\item $v$ runs over all the primes of $K$.
\end{itemize}

Poitou (1966, 1967) proved similar theorems for finite Galois modules at about the same
time as Tate, and so the duality theorems are usually credited to both. Serre
alerted each of Poitou and Tate to the work of the other, but they do not seem
to have had any direct contact.\footnote{``Je voulais te signaler que Poitou a
travaill\'e \`a peu pr\`es dans la m\^eme direction que toi\ldots J'esp\`ere
qu'au moins l'un de vous deux r\'edigera quelque chose!'' Serre, letter to Tate,
21 June 1963.}

\subsection{Tate's 1963 proof}

We seem not to know Tate's original proof of his global duality theorems, but
in a letter\footnote{The letters of Serre and Tate can be found in Volume I of 
Correspondance Serre-Tate 2015 (see References).} to Serre (25 April 1963), he observed that the nine-term sequence
can be obtained as an Ext-sequence. Specifically, on applying $\Ext(M,-)$ to
the short exact sequence
\[
0\rightarrow L^{\times}\rightarrow(\text{id\`{e}les $J_L$ of }L)\rightarrow
(\text{id\`{e}le classes $I_L$ of }L)\rightarrow0,
\]
and passing to the direct limit over the finite extensions $L$ of $K$ in
$\bar{K}$, we obtain an exact sequence
\[
\begin{tikzcd}[ampersand replacement=\&, column sep=scriptsize]
0\arrow{r}\&\Ext_{K}^{0}(M,\mathbb{G}_{m})\arrow{r}\&\Ext_{K}^{0}(M,J)\arrow{r}\&\Ext_{K}^{0}(M,C)\arrow{d}\\
\&\Ext_{K}^{1}(M,C)\arrow{d}\&\Ext_{K}^{1}(M,J)\arrow{l}\&\Ext_{K}^{1}(M,\mathbb{G}_m)\arrow{l}\\
\&\Ext_{K}^{2}(M,\mathbb{G}_{m})\arrow{r}\&\Ext_{K}^{2}(M,J)\arrow{r}
\&\Ext_{K}^{2}(M,C)\arrow{r}\&0
\end{tikzcd}
\]
that can be identified with the previous nine-term sequence by switching $M$
and $M^{D}$ and modifying the groups at the infinite primes.

Tate gave a detailed account of this proof in a letter to Tonny Springer (13
January 1966), which he intended to publish in the ``Book of the Brighton
conference on class field theory'', but which was not included. However, the
letter was widely distributed and eventually published in his Collected Works
(Part I, p.~679).

\subsection{Global duality (variant)}

We state a variant of Tate's global duality theorem in which the products over
all primes of $K$ are replaced by a direct sum over a finite set $S$ of
primes. The previous version can be obtained from this version by passing to a
direct limit over the sets $S$.

There is an exact sequence
\[
\begin{tikzcd}[column sep=scriptsize]
0 \arrow{r}& H^{0}(K_S,M) \arrow{r}{\beta^{0}}&
\bigoplus\limits_{v\in S}H^{0}(K_v,M)\arrow{r}{\gamma^{0}}
&H^{2}(K_S,M^{D})^{\ast}\arrow{d}\\
&H^{1}(K_v,M^{D})^{\ast}\arrow{d}
&\bigoplus\limits_{v\in S} H^{1}(K_v,M)\arrow{l}[swap]{\gamma^{1}}
&H^{1}(K_S,M)\arrow{l}[swap]{\beta^{1}}\\
&H^{2}(K_S,M) \arrow{r}{\beta^{2}}
&\bigoplus\limits_{v\in S}H^{2}(K_v,M)\arrow{r}{\gamma^{2}}
&H^{0}(K_S,M^{D})^{\ast}\arrow{r}
&0,\end{tikzcd}
\]
where

\begin{itemize}
\item $S$ is a finite set of primes (including any archimedean primes) of the
global field $K$;

\item $K_{S}$ = maximal extension of $K$ ramified only in $S$;

\item $M$ a finite $G_{S}$-module, $G_{S}=\Gal(K_{S}/K)$;

\item $[M]$ is not divisible by the residue characteristic at any $v\notin S$;

\item $H^{r}(K_{S},M)=H^{r}(\Gal(K_{S}/K,M)$.
\end{itemize}

\bigskip We now sketch a geometric derivation of Tate's nine-term sequence in
the function field case.

\subsection{\'{E}tale duality over a curve}

For a smooth complete curve $X$ over a field $k$, we have the following dualities.

(a) When $k=\mathbb{C}$, $X(\mathbb{C})$ is a $2$-dimensional manifold, so
there is a $2$-dimensional Poincar\'e duality theorem. When $k$ is an
arbitrary algebraically closed field, we still have a $2$-dimensional duality
theorem, but now in \'etale cohomology, provided we work with finite sheaves
prime to the characteristic of $k$.

(b) When $k$ is a finite field, there is an obvious 1-dimensional duality
theorem for finite Galois modules.

(c) When $X$ is a smooth complete curve over finite field, the two dualities
add to give a 3-dimensional duality theorem.\medskip

In more detail, let $X$ be a complete smooth curve over a field $k$, and $F$ a
constructible locally free sheaf of $\mathbb{Z}{}/m\mathbb{Z}{}$-modules, $m$
not divisible by $p$ if $\mathrm{char}(k)=p\neq0$. Let $F^{\vee}%
=\mathcal{H}om(F,\mathbb{G}_{m})$.

\begin{enumerate}
\item $k$ algebraically closed. The pairing $F^{\vee}\times F\rightarrow
\mathbb{G}_{m}$ defines a duality of finite groups%
\[
H^{2-r}(X,F^{\vee})\times H^{r}(X,F)\rightarrow H^{2}(X,\mathbb{G}_{m}%
)\simeq\mathbb{Q}{}/\mathbb{Z}{}.
\]

\item $k$ a finite field, $M$ a $\Gal(\bar{k}/k)$-module, $pM=M$. Let
$M^{\vee}=\Hom(M,\mathbb{Q}{}/\mathbb{Z}{})$. The pairing $M^{\vee}\times
M\rightarrow\bar{k}^{\times}$ defines a duality of finite groups%
\[
H^{1-r}(k,M^{\vee})\times H^{r}(k,M)\rightarrow H^{1}(k,\mathbb{Q}%
{}/\mathbb{Z}{})\simeq\mathbb{Q}{}/\mathbb{Z}{}.
\]

\item $X,k,F$ as in (a), but with $k$ finite. The pairing $F^{\vee}\times
F\rightarrow\mathbb{G}_{m}$ defines a duality of finite groups%
\[
H^{3-r}(X,F^{\vee})\times H^{r}(X,F)\rightarrow H^{3}(X,\mathbb{G}_{m}%
)\simeq\mathbb{Q}{}/\mathbb{Z}{}.
\]

\end{enumerate}

\subsection{Etale duality $\leftrightarrow$ Tate duality (char $p\neq0$).}

Now consider a smooth open curve $U$ over a finite field. When we write the
exact sequence relating the usual cohomology of $U$ to its cohomology with
compact support, and replace the latter with the group given by the duality
theorem, we obtain Tate's nine-term exact sequence.

This gives a geometric explanation for the sequence, as well as a second
proof.\medskip

In more detail, let $X$, $k$, $F$ be as in (c), and let $j\colon
U\hookrightarrow X$ be an open subscheme of $X$. We get the top row of the
following diagram with $H_{c}^{r}(U,F)=H^{r}(X,j_{!}F)$, $S=X\smallsetminus
U,$ and $K_{(v)}=$field of fractions of the henselization of $\mathcal{O}%
{}_{X,v}$. As the diagram illustrates, this essentially becomes Tate's 9-term
sequence when we replace $H_{c}^{r}(U,F)$ with $H^{3-r}(U,F^{\vee})^{\ast}$
\[
\begin{tikzcd}[column sep=small,row sep=small]
\cdots\arrow{r} & H_{c}^{r}(U,F)\arrow{r}\arrow[equals]{d}
&H^{r}(U,F)\arrow{r}\arrow[equals]{dd}
&\dstyle\bigoplus_{v\in S}H^{r}(K_{(v)},F)\arrow{r}\arrow[equals]{dd}&\cdots\\
& H^{3-r}(U,F^{\vee})^{\ast}\arrow[equals]{d}\\
\cdots\arrow{r} &H^{3-r}(K_{S},M^{D})^{\ast}\arrow{r}
&H^{r}(K_{S},M)\arrow{r}
&\dstyle\bigoplus_{v\in S}H^{r}(K_{v},M)\arrow{r}&\cdots
\end{tikzcd}
\]
$M\leftrightarrow F$ on $U,\quad M=F(U),\quad G_{S}=\pi_{1}^{\text{\'{e}t}%
}(U)$.

\subsection{Artin--Verdier duality (1964)}

Below, is the theorem as Artin and Verdier originally stated it. This gives a
geometric explanation for Tate's nine-term sequence in the number field case,
as well as a second proof (but one much more difficult than Tate's proof).

\begin{theorem}
[Artin--Verdier 1964]Let $K$ be an algebraic number field, and $j\colon
U\hookrightarrow\Spec(\mathcal{O}{}_{K})$ an open subset. For any
constructible sheaf $F$ on $U$, the Yoneda pairing%
\[
\Ext_{U}^{r}(F,\mathbb{G}_{m})\times H_{c}^{3-r}(U,F)\rightarrow H_{c}%
^{3}(U,\mathbb{G}_{m})\simeq\mathbb{Q}{}/\mathbb{Z}
\]
is a nondegenerate pairing of finite groups, except possibly on the
$2$-torsion when $K$ has a real prime. Here $H_{c}^{r}(U,F)\overset{\df}{=}%
H^{r}(X,j_{!}F)$.
\end{theorem}

Note that there is no condition on finite primes. It is possible to modify the
groups $H_{c}^{i}$ so that the theorem also holds for $2$, and then deduce
Tate's global duality as in the function field case.

Artin and Verdier stated their theorem at the famous AMS conference at Woods
Hole in 1964, but did not publish the proof. Following the conference,
Mazur gave a seminar on the theorem (cf.\ Mazur 1973), and there is a detailed proof in Milne 2006
(II, \S 3).

\section{Applications to abelian varieties}

\subsection{The pairing on the Tate--Shafarevich group}

The group $A(K)$ of rational points on an abelian variety $A$ over a global
field $K$ is finitely generated. It is easy to find the torsion subgroup of
$A(K)$ (at least for elliptic curves) so the problem of computing $A(K)$ comes
down to finding a set of generators for $A(K)$ modulo torsion. By computing
rational points, one obtains a lower bound on the rank of $A(K)$. On the other
hand, the Selmer group gives an upper bound. Roughly speaking, the difference
between the bounds is measured by the Tate--Shafarevich group.

In more detail, the Tate--Shafarevich group of an abelian variety over a
global field $K$ is defined by the exact sequence%
\[
0\rightarrow\TS(A/K)\rightarrow H^{1}(K,A)\rightarrow\bigoplus_{\text{all }%
v}H^{1}(K_{v},A)).
\]
For $m\geq1$, there is an exact sequence%
\[
0\rightarrow A(K)/mA(K)\rightarrow S^{(m)}(A/K)\rightarrow\TS(A/K)_{m}%
\rightarrow0.
\]
Here $S^{(m)}(A/K)$ is the Selmer group, which provides a finite computable
upper bound on the quotient $A(K)/mA(K)$. In the early 1950s, with the help of
an electronic computer, Selmer studied the map $E(\mathbb{Q}{})/mE(\mathbb{Q}%
{})\rightarrow S^{(m)}(E/\mathbb{Q}{})$ and found empirically that for
numerous elliptic curves $E/\mathbb{Q}{}$, the difference between the
estimates on the rank of $E(\mathbb{Q})$ coming from $S^{(m)}$ and
$S^{(m^{2})}$ is even. He conjectured that this is always true.

Cassels interpreted Selmer's conjecture as saying that the order of the
Tate--Shafarevich group is a square, and conjectured that this is explained by the
existence of a canonical bilinear form%
\[
l\colon\TS(E/K)\times\TS(E/K)\rightarrow\mathbb{Q}{}/\mathbb{Z}{}%
\]
that is alternating (i.e., $l(x,x)=0$ for all $x)$ and has kernel exactly the
group of divisible elements in $\TS(E/K)$. In a series of articles, beginning
with Cassels 1959 and culminating with Cassels 1962, he proved his conjecture.

Influenced by the work of Selmer and Cassels in the 1950s, Tate 
conjectured the duality on the
Tate--Shafarevich groups for all abelian varieties, 
and he proved this in 1962.

\begin{theorem}
[Tate, Cassels for elliptic curves]Let $A$ be an abelian variety over a number
field $K$, and let $A^{\prime}$ be its dual. There is a canonical bilinear
pairing%
\[
l\colon\TS(A/K)\times\TS(A^{\prime}/K)\rightarrow\mathbb{Q}/\mathbb{Z}%
\]
whose kernels are exactly the divisible parts of the groups. If $E$ is
a $K$-rational divisor on $A$ and $\varphi\colon
A\rightarrow A^{\prime}$ is the homomorphism $a\mapsto \mathrm{Cl}(E_a-E)$ it defines, then%
\[
l(x,\varphi(x))=0\text{ for all }x.
\]

\end{theorem}

In particular, the Tate--Shafarevich group of the Jacobian of
a curve $C$ over a number field $K$ has order a square if $C$ has a
$K$-rational point, but not in general 
otherwise.\footnote{See Poonen and Stoll 1999 for examples.}

Tate stated his theorem in his talk at the 1962 ICM, and sketched the proof in
a letter to Serre (28 July 1962). There is a detailed proof in Milne 1986, I,
\S 6. The proofs of Cassels and Tate apply also to global fields $K$ of
characteristic $p\neq0$ provided one ignores the $p$-components of the groups.

The pairing $l$ has become known as the Cassels--Tate pairing. 

\subsection{The Birch--Swinnerton-Dyer (BSD) conjecture}

\hfill\begin{minipage}{4.0in}
\textit{Birch and Swinnerton-Dyer made their conjecture for
elliptic curves \ldots over the rational numbers.
It seemed to me that the natural setting for them is for abelian varieties of any dimension, defined over any global
field.}\footnotemark\\
\hspace*{\fill}Tate, Collected Works, Part I, p.~237.
\end{minipage}
\footnotetext{The referee adds: Swinnerton-Dyer was adamant that he and Birch
 were only making the conjecture for elliptic curves over $\mathbb{Q}$ with CM; 
others could conjecture what they wanted.}

\bigskip For an elliptic curve $E$ over $\mathbb{Q}$, Birch and
Swinnerton-Dyer conjectured that the $L$-series $L(E,s)$ has a zero of order
the rank $g$ of $E(\mathbb{Q})$ at $s=1$ and, when the rank is zero, $L(E,0)$ is equal
to an expression involving the order of the Tate--Shafarevich group $\TS(E)$
of $E$. But what if $g>0$? Birch, in the proceedings of
a 1963 conference, wrote,

\begin{quote}
Tate has given a fairly detailed conjecture. One feels that $L(E,s)/(s-1)^{g}$
at $s = 1$ should give a measure of the density of rational points on the
curve $E$; so first one must decide how to measure this density. To do this,
one needs a canonical measure for the size of the generators of $E(\mathbb{Q}%
)$. This has been provided; I can give no reference beyond a letter from Tate
to Cassels. (Birch 1965).
\end{quote}

In the letter, Tate deduced, using only standard properties of heights, a
remarkably short proof of a conjecture of N\'eron that there exists a canonical quadratic 
height on abelian varieties.

\begin{conjecture}
[BSD, Tate]Let $A$ be an abelian variety over a global field $K$. Then%
\[
\lim_{s\rightarrow1}\frac{L^{\ast}(A,s)}{(s-1)^{\mathrm{rk}(A(K))}}%
=\frac{[\TS(A/K)]\cdot D}{[A(K)_{\mathrm{tors}}][A^{\prime}(K)_{\mathrm{tors}}%
]},
\]
where $A^{\prime}$ is the dual abelian variety and $D$ is the discriminant of
the N\'{e}ron-Tate height pairing $A(K)\times A^{\prime}(K)\rightarrow
\mathbb{R}{}$.
\end{conjecture}

\begin{aside}
While Tate was confident of the conjecture, not everyone was. Indeed, it was a
leap to take a statement based on calculations concerning CM elliptic curves over
$\mathbb{Q}$ and extend it to all abelian varieties over global fields, including $p$-phenomena in characteristic $p$. In 1967, Tate received a letter
from Andr\'e Weil claiming an example of an elliptic curve over a global
function field with infinite Tate--Shafarevich group,\footnote{My
recollection. It would be interesting to know if the letter still exists.} but
by then I had already proved that the group was finite in the case considered
by Weil.
\end{aside}

One important application of the duality theorems is the isogeny invariance of
the BSD conjecture.

\begin{theorem}
[Tate, Cassels for elliptic curves]Let $A$ and $B$ be abelian varieties over a
global field. If $A$ and $B$ are isogenous by an isogeny of degree prime to
the characteristic, then BSD is true for both if it true for one.
\end{theorem}

\noindent The proof uses  (for $S$ a suitable finite set of primes of $K$)

\begin{itemize}
\item Tate's global duality theorem for
\[
M\overset{\df}{=}\Ker(A(K_{S})\rightarrow B(K_{S})).
\]

\item Cassels--Tate duality (for $A$ and $B$)
\[
\TS(K,A)\times\TS(K,A^{\prime})\rightarrow\mathbb{Q}/\mathbb{Z}%
\]

\item Euler--Poincar\'{e} formula (Tate),%
\[
\frac{\lbrack H^{0}(K_{S},M)][H^{2}(K_{S},M)]}{[H^{1}(K_{S},M)]}%
=\prod_{v\text{\textrm{ arch}}}\frac{[H^{0}(K_{v},M)]}{|[M]|_{v}}\text{.}%
\]

\end{itemize}

More precisely, using the first two assertions, one finds that BSD for $A$ and
$B$ are equivalent if and only if the third assertion is true, so Tate proved
it (not without difficulty).

In the summer of 1967, I asked Tate how to prove his theorem. My recollection
is that he was able to write down a complete proof without looking anything
up, and I included the proof in my book (Milne 1986).

\begin{aside}
Tate's local and global duality theorems may be his most cited results, but 
it would be difficult to document this because there is no standard source
for them.  
As but one application, I mention that Wiles's article (1995) on Fermat's Last Theorem, makes 
crucial use of Poitou--Tate duality and Tate's Euler--Poincar\'e formula. Dick Gross reminded me
that Tate's local and global duality theorems also figure crucially in
Kolyvagin's work on elliptic modular curves and in subsequent work on Euler systems. 
He writes:
\begin{quote}
John and I read Kolyvagin's original paper in the Russian. 
We were told to find the mistake, but eventually concluded that it was brilliant. 
When I said to John that it was mostly an application of his duality theorems, 
he replied ``Well, I'm glad that someone found a use for it''.
\end{quote}
\end{aside}

\section{Local flat duality}

Tate worked with Galois cohomology, which is inadequate for the study of
$p$-phenomena in characteristic $p$. He largely left the study at $p$ to his students.

To illustrate the difference between Galois (= \'{e}tale) cohomology and flat
cohomology, consider a commutative group scheme $G$ over a field $K$, and let
$L$ be a finite extension of $K$. From the system
\[
\begin{tikzcd}[ampersand replacement=\&]
L\arrow[r, shift left=0.4ex,"a\mapsto 1\otimes a"]
\arrow[r,shift right=0.4ex,"a\mapsto a\otimes 1"']
\&L\otimes_{K}L\arrow[r, shift left=0.6ex]\arrow[r]\arrow[r,shift right=0.6ex]
\&L\otimes_{K}L\otimes_{K}
L\arrow[r, shift left=0.9ex]\arrow[r, shift left=0.3ex]\arrow[r, shift right=0.3ex]
\arrow[r, shift right=0.9ex]
\&\cdots
\end{tikzcd}
\]
we get a complex,
\[
\begin{tikzcd}[ampersand replacement=\&,column sep=small]
G(L)\arrow{r}\&G(L\otimes_K L)\arrow{r}\&G(L\otimes_K L\otimes_K L)\arrow{r}\&\cdots
\end{tikzcd}
\]
whose $r$th cohomology group we denote $H^{r}(L/K,G)$. Then%
\begin{align*}
H_{\text{\'{e}t}}^{r}(K,G)  &  =\varinjlim_{L\subset\bar{K},\,L\text{
separable over }K}H^{r}(L/K,G)\\
H_{\text{fl}}^{r}(K,G)  &  =\varinjlim_{L\subset\bar{K}}H^{r}(L/K,G).
\end{align*}
When $L/K$ is separable, $L\otimes\cdots\otimes L$ is a product of fields, and
so $G(L\otimes\cdots\otimes L)=0$ if $G$ is infinitesimal. On the other hand,
if $L/K$ is inseparable, $L\otimes_{K}L$ may have nilpotents. When $G$ is
smooth, the two groups coincide,%
\[
H_{\text{\'{e}t}}^{r}(K,G)=H_{\text{fl}}^{r}(K,G).
\]

\subsection{Finite coefficients}

A student of Tate, Steve Shatz, took up the problem of extending Tate's local
duality to local fields of characteristic $p\neq0$. He succeeded in proving a
flat duality theorem for finite group schemes in 1962, but the corresponding
theorem for abelian varieties was not proved until almost 10 years later.

\begin{theorem}
[Shatz thesis, 1962]Let $K$ be local field of characteristic $p$ (finite
residue field). Let $N$ be a finite commutative group scheme over $K$, with
Cartier dual $N^{D}$. For all $r$, the cup-product pairing%
\[
H^{r}(K,N)\times H^{2-r}(K,N^D)\rightarrow H^{2}(K,\mathbb{G}_{m}%
)\simeq\mathbb{G}_{m}%
\]
is a perfect duality\footnote{By this, I mean that the pairing realizes each
group as the Pontryagin dual of the other.} of locally compact groups.
\end{theorem}

The theorem was published as Shatz 1964, and there is also a good explanation
of it in Chapter 6 of Shatz 1972. 

\subsection{Abelian varieties}

In the 1940s, Andr\'e Weil developed a robust theory of algebraic varieties,
including abelian varieties, over arbitrary fields. This theory had difficulty
handling $p$-phenomena in characteristic $p$, essentially because it did not
allow nilpotents. For example, in the algebraic geometry of that period, there
are many maps that should be isomorphisms, but are only proved to be purely
inseparable.\bigskip

Cartier (1960) and Nishi (1959) independently extended Weil's theory of
abelian varieties to cover $p$-phenomena in characteristic $p$. Let
$\alpha\colon A\rightarrow B$ be an isogeny of abelian varieties over a field
$K$ and $\alpha^{\prime}\colon B^{\prime}\rightarrow A^{\prime}$ the dual
isogeny. In the exact sequences%
\begin{align*}
0 \rightarrow N\longrightarrow &  A\overset{\alpha}{\longrightarrow
}B\rightarrow0\\
0 \rightarrow N^{D}\longrightarrow &  B^{\prime}\overset{\alpha^{\prime
}}{\longrightarrow}A^{\prime}\rightarrow0,
\end{align*}
the finite group scheme $N^{D}$ is the Cartier dual of $N$,%
\[
N^{D}\overset{\df}{=}\mathcal{H}{}om(N,\mathbb{G}_{m}).
\]
Moreover, the canonical map $A\rightarrow A^{\prime\prime}$ from $A$ into its
double dual is an isomorphism, and the second sequence can be obtained as the
$\mathcal{E}{}xt(-,\mathbb{G}_{m})$ sequence of the first.

\begin{theorem}
[Milne 1970/1972]Let $A$ be an abelian variety over a local field $K$ and
$A^{\prime}$ the dual abelian variety. For all $r$, Tate's pairing
\[
H^{r}(K,A)\times H^{1-r}(K,A^{\prime})\rightarrow H^{2}(K,\mathbb{G}%
_{m})\simeq\mathbb{Q}/\mathbb{Z}{}%
\]
is a perfect duality (of locally compact groups).
\end{theorem}

The proof is based on Shatz's theorem. It uses the theory of N\'{e}ron minimal
models to pass to the case that $A$ and $A^{\prime}$ have semistable reduction.

Note that the statement of the theorem is exactly the same as that of Tate's
theorem --- in particular, the groups are Galois cohomology groups --- except
that it holds also for the $p$ components of the groups in characteristic $p$.

\section{Global flat duality}

As we have seen, the Poitou-Tate duality theorems can be interpreted as
duality theorems in the \'{e}tale cohomology of rings of integers in number
fields or curves over finite fields. Before long, mathematicians found that
they needed more general results, treating, for example, finite groups schemes
whose order is divisible by some residue characteristics. These require the
use of the flat topology.\footnote{By the flat topology, I mean the fppf
topology.} The author needed such a theorem in the curve case in his 1967 thesis,
and Mazur (1972) needed such a theorem in the number field case for his study
of the rational points of abelian varieties in towers of number fields.

Such theorems have been widely used. Rather than attempting to untangle their
history, I simply provide the most general result, as stated 
and proved in a recent article, Demarche and Harari 2019.

\begin{theorem}
[Artin--Mazur--Milne duality]Let $U$ be a nonempty open subset of either (a) the
spectrum of the ring of integers in a number field or (b) a complete smooth
geometrically connected curve over a finite field, and let $N$ be a finite
flat commutative group scheme over $U$. With a suitable definition of flat
cohomology with compact support, the canonical pairing%
\[
H_{c}^{3-r}(U,N)\times H^{r}(U,N^{D})\rightarrow H_{c}^{3}(U,\mathbb{G}%
_{m})\simeq\mathbb{Q}{}/\mathbb{Z},\quad0\leq r\leq3,
\]
is a perfect duality between the profinite group $H_{c}^{3-r}(U,N)$ and the
discrete group $H^{r}(U,N^{D})$ (the groups are finite in the number field case).
\end{theorem}

Using
the theorem, it is possible to extend some of the earlier theorems on abelian
varieties to the $p$ part in characteristic $p$.

\section{Interlude: arithmetic geometry in the 1960s}

The recognition in the late 1950s that algebraic geometry was the study of
schemes, and the vigorous development of scheme theory by Grothendieck and his
collaborators, enabled the great reformulation of 
arithmetic geometry that took place in the 1960s. It was during the 1960s that the foundations were laid for
the proof of the Weil conjectures and for Grothendieck's theory of motives. It
was also during the 1960s that Langlands gave birth to his theory of automorphic
representations and what we now call the Langlands program. 
Shimura varieties provided a common interest for the arithmetic geometers and the analysts.

\bigskip The 1960s was also something of a golden age for the Harvard
mathematics department. There was much cooperation between Cambridge and
Paris: both Grothendieck and Serre visited the department for substantial
periods in the 1960s, and Harvard mathematicians were frequent visitors to Paris.

Zariski's students Artin, Hironaka, and Mumford completed their degrees in
1960 and 1961. In 1962 alone, there were seminars by Hironaka (resolution of
singularities), Mumford (moduli of polarized abelian varieties), Artin (in
which \'{e}tale cohomology went from being an idea of Grothendieck to a
mathematical theory), Grothendieck (Pic, local cohomology), Tate, Kodaira,
Thompson, \ldots.\footnote{Tate, letter to Serre. April 1962} Meanwhile, in
Paris, N\'{e}ron was explaining his new theory of integral models of abelian varieties.

In the summer of 1964, there was the famous month-long conference on algebraic
geometry at Woods Hole, organized by Zariski, and attended by all the major
figures in the field except Grothendieck. This is where Tate explained his
conjectures, Artin and Verdier stated their duality theorem, Serre and Tate
stated their lifting theorem, \ldots

\bigskip I was a student at Harvard 1964--67. At the time, Brauer and Zariski
were still active, and there were also Hironaka, Mazur, Mumford, and Tate.
Tate spent the academic year 1965/66 in Paris, during which he wrote his
article with Shafarevich, proved an important case of the Tate conjecture,
and, as I shall describe, gave a Bourbaki seminar.

When he returned in the summer of 1966, I told him that I had been studying
flat cohomology and he suggested that I prove that the Tate--Shafarevich
group\footnote{In my presence, Tate always called it the Shafarevich group,
while I stubbonly stuck to Tate--Shafarevich group, until one day we both
switched to \textquotedblleft Shah\textquotedblright. Peace reigned.} is
finite. In 1966 the group was hidden in a dense fog which has scarcely
lifted, so Tate's suggestion requires explanation. This I provide in the next section.

\section{The Artin--Tate conjecture}

In this section, I report on Tate's Bourbaki Seminar of February 1966.

\subsection{Tate conjecture for surfaces over finite fields}

Let $X$ be a smooth complete surface over $\mathbb{F}_{q}$. It follows from
the Lefschetz trace formula in \'{e}tale cohomology that%
\[
\zeta(X,s)=\frac{P_{1}(X,q^{-s})P_{3}(X,q^{-s})}{(1-q^{-s})P_{2}%
(X,q^{-s})(1-q^{2-s})},\quad P_{i}(X,T)\in\mathbb{Z}{}[T].
\]

\begin{conjecture}
[Tate]The order of the pole of $\zeta(X,s)$ at $s=1$ is the rank of the
N\'{e}ron--Severi group $\mathrm{NS}(X)$ of $X$.
\end{conjecture}

Note that the order of the pole of $\zeta(X,s)$ at $s=1$ is equal to the order
of zero of $P_{2}(X,q^{-s})$ at $s=1$.

In a letter to Serre (11 June 1963), Tate said that the conjecture should be
formulated for schemes of finite type over $\mathbb{Z}{}$, and

\begin{quote}
\ldots\ \textit{most important} it should get a refinement relating the
highest coefficient of the principal part of $\zeta$ at the pole to a
discriminant attached to the group of N\'{e}ron--Severi type whose rank is the
order of the pole and to the order of a Shafarevich or Brauer-type group, just
as Birch and Swinnerton-Dyer are attempting to do in their special case.
\end{quote}

So what is the refined Tate conjecture for smooth complete surface $X$ over a
finite field $k$? We give Tate's answer in the next subsection.

\subsection{The Artin--Tate conjecture}

In collaboration with Mike Artin, Tate (1966) studied the question by mapping $X$ to
a curve $C$ in such a way that the generic fibre $X_{\eta}\rightarrow\eta$ is
smooth. Hence, $X_{\eta}$ is a smooth curve over the global function field
$K\overset{\df}{=}k(C)$, and their idea was to investigate what the BSD
conjecture for the Jacobian of $X_{\eta}$ said about $X$.

\begin{minipage}{2.0in}
\begin{tikzcd}[column sep=small]
X\arrow{d}{f}
&X_{\eta}\arrow{d}{\renewcommand{\arraystretch}{0.7} \begin{array}
[c]{l}\text{generic}\\
\text{fibre}\end{array}}\arrow{l}
&J=\mathrm{Jac}(X_\eta)\arrow[-]{d}\\
C&\eta\arrow{l}& K=k(C)
\end{tikzcd}
\end{minipage}\hspace{1cm} \begin{minipage}{2.5in}
Base field $k=\mathbb{F}_q$ (finite)\\
$X$ smooth projective surface\\
$C$ smooth projective curve\\
$f$ has smooth generic fibre $X_\eta/K$.
\end{minipage}

\medskip\noindent The result is summarized in the next two statements.

\begin{conjecture}
[BSD]%
\[
\lim_{s\rightarrow1}\frac{L^{\ast}(J,s)}{(s-1)^{\mathrm{rk}(J(K))}}%
=\frac{[\TS(J)]\cdot D}{[J(K)_{\mathrm{tors}}]^{2}},
\]

\end{conjecture}

\begin{conjecture}
[Artin--Tate]%
\[
\lim_{s\rightarrow1}\frac{P_{2}(X,q^{-s})}{(1-q^{1-s})^{\mathrm{rk}%
(\mathrm{NS}(X)))}}=\frac{\lbrack\Br(X)]\cdot D}{q^{\alpha(X)}[\mathrm{NS}%
(X)_{\mathrm{tors}}]^{2}},
\]
Here $\Br(X)$ is the Brauer group of $X$, $D$ is the discriminant of the
intersection pairing on $\mathrm{NS}(X)$, and $\alpha(X)=%
\chi(X,\mathcal{O}_{X})-1+\dim\Pic^{0}(X)$.
\end{conjecture}

The terms in the two conjectures roughly correspond. For example, using
the Leray spectral sequence for $f$, Artin proved a relation between
 $\Br(X)$ and $\TS(J)$, from which it follows that, for any prime $\ell\neq p$, 
the $\ell$-primary components of the two groups are either both finite or
both infinite (cf.\ Grothendieck 1968, \S 4).

\begin{conjecture}
[\textnf{d}]In the situation of the diagram, the two conjectures are equivalent.
\end{conjecture}

In his Bourbaki seminar, Tate stated four conjectures: (A) is the first form
of BSD for abelian varieties over global fields and (B) the full form; (C) is
what we now call the Artin--Tate conjecture, and (d) is the conjecture that,
in the context of the above diagram, Conjectures (B) and (C) are equivalent.
The last conjecture gets only a small ``d'' because, rather than being a deep
conjecture, it is a conjectural relation between deep conjectures.

\subsection{The theorems of Artin and Tate}

Tate's Bourbaki talk contained more than conjectures. He also proved the
following theorems (joint with Artin).

Let $X$ be a smooth complete surface over $\mathbb{F}_{q}$, $q=p^{a}$.

\begin{theorem}
[5.1]Let $\ell\neq p$. There is a canonical skew-symmetric form%
\[
\Br(X)(\ell)\times\Br(X)(\ell)\rightarrow\mathbb{Q}/\mathbb{Z}%
\]
on the $\ell$-primary component of $\Br(X)$ 
whose kernel consists exactly of the divisible elements.
\end{theorem}

\begin{theorem}
[5.2]Let $\ell\neq p$. The group $\Br(X)(\ell)$ is finite if and only if the
Tate conjecture holds for $X$, in which case the order of 
$\Br(X)(\ell)$ is as predicted by the
Artin--Tate conjecture.
\end{theorem}

Tate concluded his talk with the statement.

\begin{quote}
The problem of proving analogs of theorems 5.1 and 5.2 for $\ell=p$ should
furnish a good test for any $p$-adic cohomology theory, and might well serve
as a guide for sorting out and unifying the various constructions that have
been suggested: Serre's Witt vectors, Dwork's Banach spaces, the raisings via
special affines of Washnitzer-Monsky, and Grothendieck's flat cohomology for
$\mu_{p^{n}}$.
\end{quote}

Indeed, by the time we were able to prove the $p$-analogs of theorems 5.1 and 5.2, we
did know the ``correct'' $p$-adic cohomology theories. In the rest of the
article, I shall explain this and also how Conjecture (d) was proved.\label{promise}

\subsection{The case $\ell=p$ (product of two curves)}

When Tate arrived back at Harvard, not long after giving his Bourbaki talk,
and I told him that I had been studying flat cohomology, my thesis topic
presented itself: it was to understand the $p$-part of the Artin--Tate conjecture and
(a related question) the $p$-part of the BSD conjecture over a global field of characteristic $p$.

For a while I made no progress, but, at some point, Tate suggested that I look
at an example where the conjecture predicted that the Brauer group is trivial,
because it may be easier to prove that a group is trivial than to prove that
it is finite. In special cases, the Artin--Tate conjecture takes on a simple
and explicit form.

For example, when $E_{1}$ and $E_{2}$ are nonisogenous elliptic curves over
$\mathbb{F}_{q}$, the Artin--Tate conjecture says that%
\[
\lbrack\Br(E_{1}\times E_{2})]=(N_{1}-N_{2})^{2},\quad N_{i}=[E_{i}%
(\mathbb{F}{}_{q})].
\]
Note that this predicts that the order of the Brauer group is a square, as
expected. Also that, while the Brauer group may be trivial, its order cannot be
zero, and so the equation predicts that two elliptic curves over a finite
field with the same number of rational points must be isogenous. This can be
considered the zeroth case of the Tate conjecture on algebraic cycles.\footnote{This case follows
from results on the lifting of elliptic curves, proved by Deuring in the
1930s, as was noted by Mumford and also by Birch and Swinnerton-Dyer.
 See Tate 1965, p.\ 99 and also the letters from Tate to 
Serre, 9 May 1962 and 18 June 1962.}

For the case of the product of two elliptic curves, I eventually concluded
that the key was a certain flat cohomology group. More precisely,
I concluded that the key to understanding the $p$-analog in the case 
$X=E_{1}\times E_{2}$, is the flat cohomology group%
\[
H_{\mathrm{fl}}^{1}(E_{1},E_{2,p}),\quad E_{2,p}\overset{\df}{=}\Ker(E_{2}%
\overset{p}{\longrightarrow}E_{2}).
\]

When I explained this to Tate, I had no idea that anyone knew anything about
the finite group scheme $E_{p}\overset{\df}{=}\Ker(p\colon E \rightarrow E)$, but, of course, 
Tate did. When he explained its structure to me I was able, on the spot,
to prove the finiteness of $\Br(X)(p)$ in some cases.

Eventually, in my thesis (1967), I proved the $p$-analogs of the theorems 5.1
and 5.2 of Artin and Tate for the product of two curves. Since Tate had proved the Tate
conjecture in that case, this gave the following theorem.

\begin{theorem}
[Tate, Artin--Tate, Milne]The Artin--Tate conjecture holds for the product of
two curves.
\end{theorem}

At the same time, I proved that the full BSD conjecture holds for constant
abelian varieties over global fields --- in particular, that their
Tate--Shafarevich groups are finite. (An abelian variety over a global
function field is constant if it is defined by equations with coefficients in
the field of constants).\footnote{Weil's example, mentioned earlier, was an elliptic curve with
constant $j$-invariant. Thus, the curve need not be constant, but becomes
constant after a finite extension of the base field. However, if a Tate--Shafarevich group
becomes finite after a finite extension of the base field, then it was already
finite.}

These are interesting results, but not yet what I promised.

\subsection{Key step in proof of $p$ case: duality!}

Although it seems trivial now, what gave me the most trouble in my thesis was
proving a duality theorem for finite flat group schemes over a curve.

\begin{theorem}
[Artin--Milne 1976]Let $X$ be a smooth complete curve over a finite field $k$.
For a finite flat commutative group scheme $N$ over $X$ and its Cartier dual
$N^{D}$, the cup-product pairing%
\[
H^{r}(X,N)\times H^{3-r}(X,N^{D})\rightarrow H^{3}(X,\mathbb{G}_{m}%
)\simeq\mathbb{Q}{}/\mathbb{Z}{}%
\]
is a perfect duality.
\end{theorem}

For my thesis, I needed the result only for the pairs $(\alpha_{p},\alpha
_{p})$ and $(\mathbb{Z}/p\mathbb{Z},\mu_{p})$. Note that the pairing
\[
(m,\zeta)\mapsto\zeta^{m}\colon\mathbb{Z}{}/p\mathbb{Z}{}\times\mu
_{p}\rightarrow\mathbb{G}_{m}.
\]
realizes each of $\mathbb{Z}/p\mathbb{Z}$ and $\mu_{p}$ as the Cartier dual of
the other. Following is a sketch of the proof in this case.

There is an Artin--Schreier sequence, exact on $X_{\mathrm{et}}$,
\begin{equation}
\begin{tikzcd} 0\arrow{r}&\mathbb{Z}/p\mathbb{Z}\arrow{r}&\mathcal{O}_{X}\arrow{r}{x\mapsto x^{p}-x}&\mathcal{O}_{X}\arrow{r}&0, \end{tikzcd} \tag{*}%
\end{equation}
and a Kummer sequence, exact on $X_{\text{fl}}$,%
\[
\begin{tikzcd}[ampersand replacement=\&]
1\arrow{r}
\&\mu_{p}\arrow{r}\&\mathbb{G}_{m}\arrow{r}{x\mapsto x^{p}}\&\mathbb{G}_{m}\arrow{r}\&1.
\end{tikzcd}
\]
On applying the morphism of sites $f\colon X_{\text{fl}}\rightarrow X_{\text{et}}$
defined by the identity map, we deduce that
\[
R^{i}f_{\ast}\mu_{p}\simeq\left\{
\begin{array}
[c]{ll}%
\mathcal{O}{}_{X}^{\times}/\mathcal{O}{}_{X}^{\times p}\overset{\df}{=}%
\nu(1) & \text{if }i=1\\
0\quad & \text{otherwise,}%
\end{array}
\right.
\]
so
\begin{equation}
H_{\text{fl}}^{i}(X,\mu_{p})\simeq H_{\text{et}}^{i-1}(X,\nu(1)). \tag{**}%
\end{equation}
From the exact sequence ($C$ is the Cartier operator)%
\[
\begin{tikzcd}[ampersand replacement=\&,column sep=scriptsize,row sep=small]
0\arrow{r}
\&\mathcal{O}_{X}^{\times}\arrow{r}{h\mapsto h^{p}}\&\mathcal{O}_{X}^{\times}\arrow{rr}{h\mapsto dh/h}\arrow[two heads]{rd}
\&\&\Omega_{X}^{1}\arrow{r}{C-1}
\&\Omega_{X}^{1}\arrow{r}\&0\\
\&\&\& \nu(1)\arrow[hook]{ru}
\end{tikzcd}
\]
we can extract a short exact sequence%
\begin{equation}
0\rightarrow\nu(1)\rightarrow\Omega_{X}^{1}\overset{C-1}{\longrightarrow
}\Omega_{X}^{1}\rightarrow0. \tag{***}%
\end{equation}
There is a $1$-dimensional duality between the Zariski (= \'etale) cohomologies of
$\mathcal{O}{}_{X}$ and $\,\Omega_{X}^{1}$. Using (*) and (***), we deduce a
$2$-dimensional duality between the \'{e}tale cohomologies of $\mathbb{Z}%
{}/p\mathbb{Z}{}$ and $\nu(1)$. Finally, using (**), we obtain a
$3$-dimensional duality between the flat cohomologies of $\mathbb{Z}%
/p\mathbb{Z}{}$ and $\mu_{p}$.

\subsubsection{The Cartier operator}

For a smooth variety $X$ over a perfect field $k$, Cartier (1957) showed that
there is a (unique) family of maps%
\[
C\colon\Omega_{X/k,\text{closed}}^{r}\rightarrow\Omega_{X/k}^{r}%
\]
such that

\begin{itemize}
\item $C(\omega+\omega^{\prime})=C(\omega)+C(\omega^{\prime}),\quad
C(h^{p}\omega)=hC(\omega),$

\item $C(\omega\wedge\omega^{\prime})=C(\omega)\wedge C(\omega^{\prime})$,

\item $C(\omega)=0$ $\iff\omega$ is exact,

\item $C(dh/h)=dh/h.$
\end{itemize}

\medskip\noindent For curves, the Cartier operator was defined by Tate (1952)
in a paper in which he studied how the genus of a curve changes under
extension of the base field.\footnote{The reader may object that the genus
does not change under extension of the base field. This is true in 2025, but things
were different in 1952. Consider a complete normal curve $X$ over a field $k$
and an extension $k^{\prime}$ of $k$. The curve $X^{\prime}$ obtained by
extending the base field to $k^{\prime}$ does indeed have the same genus as
$X$, but it may no longer be normal, for example, its structure sheaf may
acquire nilpotents. In 1952, by the extended curve one meant the associated
normal curve, whose genus may drop (but only by a multiple of $(p-1)/2$, as
proved by Tate). Tate in fact expressed himself in terms of function fields.}

\section{Flat duality (Artin's conjecture)}

To continue the story, we need another duality theorem, this time conjectured
by Artin.

\subsection{Artin's conjecture}

In an important article, Artin (1974) used flat cohomology to study
supersingular $K3$ surfaces. This led him to conjecture a duality theorem in
the flat cohomology of surfaces over fields of characteristic $p\neq0$.

Let $\pi\colon X\rightarrow\Spec k$ be a smooth complete surface over a
perfect field $k$ of characteristic $p\neq0$.

\begin{conjecture}
[Rough form]When $k$ is algebraically closed, there is a 4-dimensional duality
for the finite part of $H_{\mathrm{fl}}^{i}(X,\mu_{p})$ and a $5$-dimensional
duality for the vector space part.
\end{conjecture}

Clearly this needs to be restated in terms of derived categories. Artin
proved\footnote{Artin did not publish his proof. The statement is proved in
Bragg and Olsson 2021.} that the functor $R^{r}\pi_{\ast}\mu_{p^{n}}$ (flat
topology) is represented by a group scheme of finite type over $k$. His
conjecture concerned only these group schemes modulo infinitesimal group schemes.

\begin{conjecture}
[Precise form]There is a canonical isomorphism%
\[
R\pi_{\ast}\mu_{p^{n}}\rightarrow R\Hom(R\pi_{\ast}\mu_{p^{n}},\mathbb{Q}%
/\mathbb{Z})[-4]
\]
in the derived category of the category of commutative group schemes over $k$
modulo infinitesimal group schemes.
\end{conjecture}

\subsection{Proof of Artin's conjecture ($n=1$)}

In the proof of the flat duality theorem for curves, we saw that we should
identify the flat cohomology of $\mu_{p}$ with the \'{e}tale cohomology of the
sheaf $\nu(1)$ shifted by $1$. This idea works more generally.

Let $\pi\colon X\rightarrow\Spec k$ be a smooth complete variety of dimension
$d$ over a perfect field $k$ of characteristic $p$. Define a sheaf on
$X_{\text{\'{e}t}}$ by%
\[
\nu(r)=\Ker(C-1\colon\Omega_{X\text{,closed}}^{r}\rightarrow\Omega_{X}^{r}).
\]

\begin{theorem}
[Milne 1976]The functor $R\pi_{\ast}\nu(r)$ is representable on perfect
schemes, and there is a canonical isomorphism%
\[
R\pi_{\ast}\nu(r)\rightarrow R\Hom(R\pi_{\ast}\nu(d-r),\mathbb{Z}%
/p\mathbb{Z}{}{})[-d]
\]
in the derived category of the category of commutative group schemes killed by
$p$ modulo infinitesimal group schemes.
\end{theorem}

When $d=2$, $r=$1, this becomes Artin's conjecture for $\mu_{p}$: we have
\[
(X_{\mathrm{fl}}\overset{f}{\longrightarrow}X_{\mathrm{et}}\overset{\pi
^{\mathrm{et}}}{\longrightarrow}\Spec k)=
(X_{\mathrm{fl}}\overset{\pi^{\mathrm{fl}}}\longrightarrow\Spec k)%
\]
and
\[
Rf_{\ast}\mu_{p}=\nu(1)[-1],
\]
so
\[
R\pi_{\ast}^{\text{fl}}\mu_{p}=R\pi_{\ast}^{\text{et}}Rf_{\ast}\mu_{p}%
=R\pi_{\ast}^{\text{et}}\nu(1)[-1].
\]

\subsection{Proof of Artin's conjecture (all $n$)}

At this point I was stumped: my proof of Artin's conjecture depends on the
sheaves of differentials, which are killed by $p$ in characteristic $p$, so
how to prove the conjecture for $\mu_{p^{n}}$?

In 1974, I shared an office with Spencer Bloch at the University of Michigan.
When I explained my problem to him he said that he had constructed objects
that were just like the sheaves of differentials, except killed by $p^{n}$,
not $p$. Indeed, he had.\footnote{Published as Bloch 1977.} This was the
famous de Rham--Witt complex, which is a projective system of complexes%
\[
W_{n}\mathcal{O}_{X}\overset{d}{\longrightarrow}W_{n}\Omega_{X}^{1}%
\overset{d}{\longrightarrow}W_{n}\Omega_{X}^{2}\rightarrow\cdots,\qquad n\geq 1,
\]
of $W_{n}(\mathcal{O}_{X})$-modules.

Bloch defined the de Rham--Witt complex in order to relate $K$-theory to
crystalline cohomology, but once he had introduced it, its importance was
apparent, and it was soon developed further by others.\footnote{Initially 
by Illusie (1979) and Illusie and Raynaud (1983); more recently by 
 Bhatt, Lurie, and Matthew (2021).} Not only does it give
a new construction of crystalline cohomology, but it adds structure to it. For
example, as mentioned earlier, Serre had studied the cohomology of the sheaf
of Witt vectors on a variety, and had correctly concluded that it gives only
part of the ``good'' cohomology. With the de Rham--Witt complex, it became
possible to say exactly which part.

Using the de Rham--Witt complex, it became possible to define sheaves $\nu
_{n}(r)$ (killed only by $p^{n}$) and prove by induction from the case $n=1$
that the canonical morphism%
\[
R\pi_{\ast}\nu_{n}(r)\rightarrow R\Hom(R\pi_{\ast}\nu_{n}(d-r),\mathbb{Q}%
{}/\mathbb{Z}{})[-d]\text{.}%
\]
is an isomorphism. When $d=2$, this becomes Artin's conjecture for $\mu
_{p^{n}}$.

\section{Conclusion}

It is now possible to provide the answers to Tate's questions promised on p.~\pageref{promise}.

\subsection{The analogs for $\ell=p$ of the theorems of Artin and Tate}

Using the sheaves $\nu_{n}(r)$, it became possible to complete the proof of
the analogs for $\ell=p$ of the Theorems 5.1 and 5.2 in Tate's 1966 Bourbaki
talk.\medskip

Let $X$ be a smooth complete surface over a finite field $k$ of characteristic
$p$.

\begin{theorem}
[Milne 1975]There is a canonical skew-symmetric form
\[
\Br(X)(p)\times\Br(X)(p)\rightarrow\mathbb{Q}{}/\mathbb{Z}%
\]
whose kernel consists exactly of the divisible elements.
\end{theorem}

\begin{theorem}
[Milne 1975, completing Artin and Tate]The following are equivalent.

\begin{enumerate}
\item The Tate conjecture holds for $X$.

\item For some prime $l$ ($l=p$ is allowed), $\Br(X)(l)$ is finite.

\item The Artin--Tate conjecture holds for $X$ (including the $p$ part).
\end{enumerate}
\end{theorem}

\begin{aside}
In his 1962 ICM talk, Tate said that he suspects that the form on the Brauer
group is not only skew-symmetric, but alternating, so that the order of the
Brauer group is a square if finite. This is true, but the proof has only
recently been completed (Carmeli and Feng 2025).
\end{aside}

\subsection{Proof of Conjecture \textnf{(d)}}

\begin{theorem}
[Kato--Trihan]Let $A$ be an abelian variety over a global function field
$K$. The following are equivalent.

\begin{enumerate}
\item The order of the zero of $L(s,A)$ at $s=1$ is the rank of $A(K).$

\item For some prime $l$, $\TS(A/K)(l)$ is finite.

\item The full BSD conjecture for $A/K$.
\end{enumerate}
\end{theorem}

For the proof, which uses global flat duality over curve,
see Kato--Trihan 2003 and Trihan and Vauclair 2024, 1.0.1.

We can now prove Tate's Conjecture (d) for the pair
\[
\begin{tikzcd}
X\arrow{d}{f}
&X_{\eta}\arrow{d}\arrow{l}
&&J=\mathrm{Jac}(X_\eta)\arrow[-]{d}\\
C&\eta\arrow{l}&& K=k(C)
\end{tikzcd}
\]
Recall that Conjecture (d) says that, in the situation of the diagram,
\[
\text{the Artin--Tate conjecture holds for $X\iff$ the full BSD conjecture
holds for $J$.}
\]
Because of the equivalences in the last two theorems, it suffices to prove
that
\[
\text{for some prime $l$, $\Br(X)(l)$ is finite $\iff$ for some prime $l$,
$\TS(J/K)(l)$ is finite.}
\]
As we noted earlier, it follows from work of Artin that, for $l\neq p$, 
the $l$-primary components of the two groups are either both finite
or both infinite.

\subsection{The \textquotedblleft good\textquotedblright\ $p$-adic cohomology
theories in characteristic $p$}

Let $X$ be a smooth complete variety over a field $k$ of characteristic
$p\neq0$. I claim that the \textquotedblleft good\textquotedblright%
\ cohomologies are%
\[
\renewcommand{\arraystretch}{1.3}\Bigg\{
\begin{array}
[c]{rl}%
\text{Weil cohomology:} & H_{\mathrm{crys}}^{r}(X/W)\simeq H^{r}(X,W\Omega
_{X}^{\bullet})\\
\text{\textquotedblleft}H_{\mathrm{fl}}^{i}(X,\mu_{p^{n}}^{\otimes
r})\text{\textquotedblright}\colon & H_{\text{\'{e}t}}^{i-r}(X,\mathbb{\nu
}_{n}(r)).
\end{array}
\]
Note that the quotation marks can be removed when $r\leq1$. The second
definition may seem too ad hoc to be convincing, but there is a second
description of it.

When we apply $R\Gamma$ to the de Rham--Witt complex of a variety, we get a
complex of $W(k)$-modules, from which we can deduce the crystalline
cohomology groups $H_{\mathrm{crys}}^{r}(X/W)$. The de Rham--Witt complex has extra structure, namely, an action
of the Raynaud ring.\footnote{In addition to the action of the Witt vectors,
the de Rham--Witt complex has operators $F, V, d$ satisfying certain
conditions. To say that an object has these operators is exactly to say that
it has an action of the Raynaud ring.} When we remember this action, the same
construction gives $\varprojlim_{n}H_{\text{\'{e}t}}^{i-r}(X,\mathbb{\nu}%
_{n}(r))$ instead $H_{\mathrm{crys}}^{r}(X/W)$.

In more detail, when we regard $R\Gamma(W\Omega_{X}^{\bullet})$ as an object
in the triangulated category with $t$-structure $\mathsf{D}^{+}(W)$,%
\[
H_{\mathrm{crys}}^{i}(X/W)\simeq\Hom_{\mathsf{D}^{+}(W)}(\1,R\Gamma
(W\Omega_{X}^{\bullet})[i]).
\]
\noindent When we regard $R\Gamma(W\Omega_{X}^{\bullet})$ as an object in the
triangulated category with $t$-structure $\mathsf{D}_{c}^{b}(R)$ ($R$ the
Raynaud ring),%
\[
\text{\textquotedblleft}\varprojlim_{n}H_{\mathrm{fl}}^{i}(X,\mu_{p^{n}%
}^{\otimes r})\text{\textquotedblright}\simeq\Hom_{\mathsf{D{}}_{c}^{b}%
(R)}(\1,R\Gamma(W\Omega_{X}^{\bullet})(r)[i])
\]
(Milne and Ramachandran 2005).

\bigskip\noindent\textsc{Postscript:} In the talk by Bhargav Bhatt following
mine at the conference, the $\nu$-sheaves re-appear as the objects
$\mathcal{O}\{r\}$ in the category of $F$-gauges. As Bhatt noted,
\v{C}esnavi\v{c}ius and Scholze (2024) used them to prove $p$-purity statements for
Picard groups and Brauer groups (conjectures of Gabber).

\bigskip
\bigskip

\begin{quote}
Tate had many Ph.D.\ students, and he took good care of us. At the conference,
Shankar Sen, who was my contemporary as a student, told how Tate once came
past his dorm out of concern for him. I can tell a similar story. At some
point, Tate decided I should finish my degree in the winter term of 1967,
which meant that there was a deadline. Specifically, I was to deliver my
completed thesis to the typist by 9\thinspace am on a certain Monday. During
the week before the deadline, I was still having trouble getting all the
pieces of my thesis to fit together. I felt so bad about this that I avoided
the mathematics department (which was then in the beautiful old building at 2
Divinity Avenue). On Saturday morning, I felt safe to resume working at my
usual place in the library in the mathematics department. To my surprise, Tate
showed up, having biked in. When I explained that I had one last statement to
prove before I could finish writing up my thesis, Tate looked at it, and said
\textquotedblleft Seems OK. You don't \textit{have} to sleep this weekend do
you?\textquotedblright, and left. In fact, I did meet the deadline, and I even
got some some sleep that weekend. I should add that if I had missed the
deadline, nothing bad would have happened --- Tate was pretty kind hearted.
\end{quote}

\clearpage
\markboth{References}{References}

\section*{References}

\begin{adjustwidth}{2em}{0em}
\setlength{\parindent}{-2em}

\hspace{-2em}
Correspondance Serre-Tate. Vol. I, Documents Math\'ematiques (Paris), 13, Soci\'et\'e{} Math\'ematique de France, Paris, 2015; MR3379329

Artin, M., 1974, Supersingular $K3$ surfaces, Ann. Sci. \'Ecole Norm.
Sup. (4) {\bf 7} (1974), 543--567; MR0371899

Artin, M. and Milne, J., 1976, Duality in the flat cohomology of curves, Invent.
Math. \textbf{35} (1976) 111--129.

Artin, M. and Verdier, J.-L., 1964, Seminar on \'{e}tale cohomology of number
fields, AMS Woods Hole Summer Institute 1964, notes available in a limited
edition only, 5pp. 

Bhatt, B., Lurie, J., and Mathew, A., 2021, 
Revisiting the de Rham--Witt complex, Ast\'erisque No. 424 (2021), viii+165 pp.; 
MR4275461

Birch, B. J., 1965, Conjectures concerning elliptic curves. Proc. Sympos. Pure Math., Vol. VIII, pp. 106-112, American Mathematical Society, Providence, RI, 1965

Bloch, S., 1977, Algebraic $K$-theory and crystalline cohomology, 
Inst. Hautes \'Etudes Sci. Publ. Math. No. 47 (1977), 187--268; MR0488288

Bragg, D., and Olsson, M., 2021, Representability of cohomology of finite 
flat abelian group schemes, arXiv:2107.11492.

Carmeli, S., and Feng, T., 2025, Prismatic Steenrod operations and arithmetic duality on Brauer groups,
arXiv:2507.13471.

Cartier, P., 1957, Une nouvelle op\'eration sur les formes diff\'erentielles, C. R. Acad. Sci. Paris {\bf 244} (1957), 426--428; MR0084497

Cartier, P., 1960, Isogenies and duality of abelian varieties, 
Ann. of Math. (2) {\bf 71} (1960), 315--351; MR0116019

Cassels, J., 1959, Arithmetic on curves of genus $1$. I. On a conjecture of Selmer, J. Reine Angew. Math. {\bf 202} (1959), 52--99; MR0109136

Cassels, J., 1962, Proof of the Hauptvermutung, J. Reine Angew. Math.
\textbf{211} (1962) 95--112.

\v Cesnavi\v cius, K., and Scholze, P., 2024, Purity for flat cohomology, Ann.
of Math. (2) {\bf 199} (2024), no.~1, 51--180; MR4681144

Demarche, C., and Harari, D., 2019, Artin--Mazur--Milne duality for 
fppf cohomology. Algebra Number Theory 13 (2019), no. 10, 2323--2357.

Illusie, L., 1979, 
 Complexe de de Rham-Witt et cohomologie cristalline. 
Ann. Sci. \'Ecole Norm. Sup. (4) 12 (1979), no. 4, 501--661. 

llusie, L.; Raynaud, M., Les suites spectrales 
associ\'ees au complexe de de Rham-Witt. 
 Inst. Hautes \'Etudes Sci. Publ. Math. No. 57 (1983), 73--212.

Kato, K., and Trihan, F., 2003, On the conjectures of Birch and
Swinnerton-Dyer in characteristic $p>0$, Invent. Math. {\bf 153}
(2003), no.~3, 537--592; MR2000469

Mazur, B., 1972, Rational points of abelian varieties with values in towers of
number fields, Invent. Math. \textbf{18}, 183--266.

Mazur, B., 1973, Notes on \'etale cohomology of number fields. Ann. Sci.
\'Ecole Norm. Sup. (4) 6 (1973), 521--552.

Milne, J., 1967, The conjectures of Birch and Swinnerton-Dyer for constant
abelian varieties over function fields, Thesis, Harvard University. 

Milne, J. 1970/72, Weil-Chatelet groups over local fields, Ann. Sci. Ecole Norm.
Sup. \textbf{3} (1970), 273--284; addendum, ibid., \textbf{5} (1972), 261--264.

Milne, J., 1975, On a conjecture of Artin and Tate, Annals of Math. 
\textbf{102} (1975) 517--533.

Milne, J., 1976, Duality in the flat cohomology of a surface, Ann. scient. ENS,
\textbf{9} (1976) 171--202.

Milne, J., 1986, Arithmetic duality theorems, Perspectives in Mathematics, 1, Academic Press, Boston, MA, 1986; MR0881804

Milne, J., 2006, Arithmetic duality theorems, second edition, 
BookSurge, Charleston, SC, 2006; MR2261462

Milne, J., and Ramachandran, N., 2005, 
The de Rham--Witt and $\mathbb{Z}_p$-cohomologies of an 
algebraic variety.  Adv. Math.  198  (2005),  no. 1, 36--42.

Nishi, M., 1959, The Frobenius theorem and the duality theorem on an abelian variety, Mem. Coll. Sci. Univ. Kyoto Ser. A. Math. {\bf 32} (1959), 333--350; MR0116020

Poitou, G. 1966, Remarques sur l'homologie des groupes profini, In \textit{Les
Tendances} \textit{G}\'{e}\textit{ometrie en Alg}\`{e}\textit{bre et
Th}\'{e}\textit{orie des Nombres}, CNRS, Paris, 201--213.

Poitou, G. 1967, \textit{Cohomologie Galoisienne des Modules Finis}, Dunod,
Paris.

Poonen, B., and Stoll, M., 1999, The Cassels--Tate pairing on polarized abelian varieties. Ann. of Math. (2) 150 (1999), no. 3, 1109--1149.

Shatz, S., 1962, Cohomology of artinian group schemes over local fields,
thesis, Harvard University.

Shatz, S., 1964, Cohomology of artinian group schemes over local fields, 
Ann. of Math. (2) {\bf 79} (1964), 411--449; MR0193093,

Shatz, S., 1972, {\it Profinite groups, arithmetic, and geometry}, Annals of
Mathematics Studies, No. 67, Princeton Univ. Press, Princeton, NJ,
1972 Univ. Tokyo Press, Tokyo, 1972; MR0347778

Tate, J., 1952, Genus change in inseparable extensions of function fields, Proc. Amer. Math. Soc. {\bf 3} (1952), 400--406; MR0047631

Tate, J., 1957, WC-groups over $\mathfrak{p}$-adic fields, S\'eminaire
Bourbaki, Expos\'e 156, D\'ecembre 1957, 13pp.

Tate, J., 1962, Duality theorems in Galois cohomology over number fields,
\textit{Proc.} \textit{Intern. Congress Math.}, Stockholm, pp234--241.

Tate, J., 1965, Algebraic cycles and poles of zeta functions, in {\it Arithmetical Algebraic Geometry (Proc. Conf. Purdue Univ., 1963)}, pp. 93--110, Harper \& Row, New York, ; MR0225778

Tate, J., 1966, On the conjectures of Birch and Swinnerton-Dyer and a
geometric analog, S\'{e}minaire Bourbaki, Expose\'e 306,   
(Reprinted in: Dix
Expos\'{e}s sur la Cohomologie des Sch\'{e}mas, North-Holland, Amsterdam, 1968).

Tate, J., Collected works of John Tate. Part I (1955-1975), 
Part II (1976-2006). Edited by Barry Mazur and Jean-Pierre Serre, 
AMS, Providence, RI, 2016.

Trihan, F., and Vauclair, D., 2024, 
A comparison theorem for semi-abelian schemes over a smooth curve, 
Mem.\ Amer.\ Math.\ Soc.\ 299 (2024), no.\b 1495.

Wiles, A., 1995, Modular elliptic curves and Fermat's last theorem. Ann.
of Math. (2) 141 (1995), no. 3, 443--551.

\end{adjustwidth}

\end{document}